\documentclass[twoside,12pt]{amsart}
\usepackage{amssymb,longtable}
\theoremstyle{definition}
\newtheorem{thm}{Theorem}[section]
\newtheorem{prop}[thm]{Proposition}

\newtheorem{lem}[thm]{Lemma}
\newtheorem{rem}[thm]{Remark}

\newtheorem{example}[thm]{Example}
\newtheorem*{acknowledgement}{Acknowledgments}

\def\rk{\mathop{\mathrm{rank}}\nolimits}

\def\Hom{\mathop{\mathrm{Hom}}\nolimits}
\title[Non-symplectic automorphisms of order 3 on $K3$ surfaces]
{Classification of non-symplectic automorphisms of order 3 on $K3$ surfaces}
\author[S.~Taki]{Shingo Taki}
\address{Graduate School of Mathematics\\ 
 Nagoya University, Chikusa-ku Nagoya 464-8602 Japan}
\email{m04022x@math.nagoya-u.ac.jp}
\address{Korea Institute for Advanced Study, Hoegiro 87, Dongdaemun-gu, Seoul 130-722, Korea}
\email{taki@kias.re.kr}
\date{\today}

\begin{document}
\bibliographystyle{amsalpha+}

\begin{abstract}
In this paper, we study non-symplectic automorphisms of order 3 on algebraic $K3$ surfaces over $\mathbb{C}$ 
which act trivially on the N\'{e}ron-Severi lattice. In particular we shall characterize their 
fixed loci in terms of the invariants of 3-elementary lattices.
\end{abstract}

\maketitle

\section{Introduction}

Let $X$ be an algebraic surface over $\mathbb{C}$. 
If its canonical line bundle $K_{X}$ is trivial and $\dim H^{1}(X, \mathcal{O}_{X})=0$ then 
$X$ is called a \textit{$K3$ surface}. 
In this paper, we study automorphisms of algebraic $K3$ surfaces.

In the following, for an algebraic $K3$ surface $X$, we denote by $S_{X}$, $T_{X}$ and $\omega _{X}$ 
the N\'{e}ron-Severi lattice, the transcendental lattice and a nowhere vanishing 
holomorphic $2$-form on $X$, respectively.

An automorphism of $X$ is \textit{symplectic} 
if it acts trivially on $\mathbb{C}\cdot \omega _{X}$.
In particular, this paper is devoted to a study of \textit{non}-symplectic automorphisms on algebraic $K3$ surfaces.  

First we recall some general results about non-symplectic automorphisms on algebraic $K3$ surfaces.
Nikulin \cite{Nikulin3} has classified non-symplectic involutions by using the classification of $2$-elementary lattices. 
In the paper, he considered fixed loci of involutions. 
And he characterizes them in terms of the invariants of 2-elementary lattices. 

More generally, we suppose that $g$ is a non-symplectic automorphism of order $I$ on $X$ such that 
$g^{\ast }\omega _{X} =\zeta _{I} \omega _{X}$ where $\zeta _{I}$ is a primitive $I$-th root of unity.
Then $g^{\ast }$ has no non-zero fixed vectors in $T_{X}\otimes \mathbb{Q}$ and hence  
$\phi (I)$ divides $\rk T_{X}$, where $\phi $ is the Euler function.  
In particular $\phi (I)\leq \rk T_{X}$ and hence $I\leq 66$ 
(\cite{Nikulin2}, Theorem 3.1 and Corollary 3.2).

The extremal case $\rk T_{X}=\phi (I)$ was determined in the result below \cite{Vorontsov}, \cite{Kondo1}. 
If $T_{X}$ is unimodular then $I\in \{66, 44, 42, 36, 28, 12\}$. 
Moreover for $I \in \{66, 44, 42, 36, 28, 12\}$, there exists a unique pair $(X_{I}, \langle g_{I} \rangle )$, 
where $X_{I}$ is a $K3$ surface and $g_{I}$ is an automorphism of order $I$.
And if $T_{X}$ is not unimodular and $g_{I}$ acts trivially on $S_{X}$ then $I\in \{3^{k}, 5^{l}, 7, 11, 13, 17, 19\}$ 
where $1\leq k\leq 3$ and $l=1, 2$.

Non-symplectic automorphisms with order 13,17 and 19 were studied by Oguiso and Zhang \cite{Oguiso-Zhang2}.
They proved that there exists a unique pair $(X_{I}, \langle g_{I} \rangle )$ of 
a $K3$ surface and an automorphisms group of order $I$ where $I \in \{13,17,19\}$. 
In particular, $(X_{I}, \langle g_{I} \rangle )$ are isomorphic to Kondo's example in \cite{Kondo1}.

Sometimes an automorphism with small order can also characterize uniquely a $K3$ surface. 
Actually non-symplectic automorphisms with order 3 were studied by Oguiso and Zhang \cite{Oguiso-Zhang1}. 
They proved that there exists a unique pair $(X, g)$ of a $K3$ surface and an automorphism of order 3 
such that the fixed locus $X^{g}=\{x\in X| g(x)=x\}$ consists of only at least 6 rational curves and some isolated  points. 

In this paper we consider non-symplectic automorphisms of order 3 acting trivially on $S_{X}$. 
In this case, $S_{X}$ is a 3-elementary lattice, that is,
$S_{X}^{\ast }/S_{X}$ is a 3-elementary group where $S_{X}^{\ast }=\Hom (S_{X},\mathbb{Z})$ 
(see Proposition \ref{vo}).

Following Nikulin's result \cite{Nikulin3}, we characterize fixed loci in terms of  the invariants of 3-elementary lattices.
The main purpose of this paper is to prove the following theorem:

\begin{thm}\label{mt}
Let $\rho$ be the Picard number of $X$  
and let $s$ be the minimal number of generators of $S_{X}^{\ast }/S_{X}$.

\begin{itemize}
\item[(1)] If $22-\rho -2s<0$, then 
$X$ has no non-symplectic automorphisms of order 3 which act trivially on $S_{X}$.

\item[(2)]  If $22-\rho -2s\geq 0$, then 
$X$ has a non-symplectic automorphism $\varphi $ of order 3 which acts trivially on $S_{X}$. 
Moreover the fixed locus $X^{\varphi }:=\{x\in X| \varphi (x)=x \}$ has the form 
\begin{equation*}
X^{\varphi }=
\begin{cases}
\{ P_{1} \} \amalg \{ P_{2} \} \amalg \{ P_{3} \} & \hspace{-2.25cm} \text{if $S_{X}=U(3)\oplus E_{6}^{\ast }(3)$} \\
\{ P_{1} \} \amalg \dots \amalg \{ P_{M} \} \amalg C^{(g)} \amalg E_{1} \amalg \dots \amalg E_{N-1} & \text{otherwise}
\end{cases}
\end{equation*}
and $M = \rho /2-1$, $g=(22-\rho -2s)/4$, $N=(6+\rho -2s)/4$, where 
we denote by $P_{i}$ an isolated point, $C^{(g)}$ a non-singular curve of genus $g$ and by $E_{j}$ a non-singular rational curve.
\end{itemize}
\end{thm}
\begin{rem}
We assume $22-\rho -2s\geq 0$.
For each $(\rho ,s)$, we can construct a $K3$ surface which has a non-symplectic automorphism of order 3 
acting trivially on $S_{X}$. (See section \ref{exa}.)

Moreover we shall give Table \ref{SandF} in the last of section \ref{form}
in which the fixed locus $X^{\varphi }$ is given for each $S_{X}$ satisfying $22-\rho -2s\geq 0$.
\end{rem}

In section \ref{class}, we shall give a classification of an even hyperbolic 3-elementary lattices admitting a primitive embedding in  $K3$ lattice. 
As a result, we get all lattices which are 
the N\'{e}ron-Severi lattice of $K3$ surfaces with non-symplectic automorphisms of order 3 which act trivially on $S_{X}$.
The main part in this paper is section \ref{form}. 
In this section, we shall prove the theorem. 
Here we use mainly the Lefschetz formula, the Hurwitz formula and the theory of elliptic surfaces due to Kodaira \cite{Kodaira}. 
In section \ref{exa}, we give examples of algebraic $K3$ surfaces with a non-symplectic automorphism of order 3 
acting trivially on $S_{X}$. 
But we do not prove the uniqueness of these examples.

\vspace{3mm}

\begin{rem}
Recently, M.~Artebani, A.~Sarti \cite{AS} independently have given a 
classification of non-symplectic automorphisms of order 3 on $K3$ surfaces 
which act trivially on the N\'{e}ron-Severi lattice by using Smith exact sequences. 
Our proof is based on geometric argument without using Smith exact sequences. 
\end{rem}

\begin{acknowledgement}
The author would like to express his gratitude to Professor Shigeyuki Kondo 
for suggesting this problem, many helpful discussions and suggestions. 
He is also grateful to Professor Michela Artebani
for informing about Example \ref{exz} and a mistake in the first version.
\end{acknowledgement}

\section{Classification of 3-elementary lattices}\label{class}
A lattice of $\rk r$ is called \textit{hyperbolic} if its signature is $(1, r-1)$.
Even hyperbolic $p$-elementary lattices were classified by \cite{R-S}. 
By using their theorem, 
we classify even hyperbolic 3-elementary lattices admitting a primitive embedding in 
an even unimodular lattice of signature $(3,19)$ satisfying $22-\rho -2s\geq 0$. 
These lattices give the N\'{e}ron-Severi lattice of $K3$ surfaces $X$ 
with non-symplectic automorphisms of order 3 which act trivially on $S_{X}$.

\vspace{3mm}

We denote by $U$ the hyperbolic lattice defined by $\begin{pmatrix} 0 & 1 \\ 1 & 0 \end{pmatrix}$ 
which is an even unimodular lattice of signature $(1,1)$, and by $A_{m}$, $D_{n}$ or $E_{l}$ an even 
negative definite lattice associated with the Dynkin diagram of type $A_{m}$, $D_{n}$ or $E_{l}$ 
($m\geq 1$,  $n\geq 4$, $l=6,7,8$). 
For a lattice $L$ and an integer $m$, $L(m)$ is the lattice whose bilinear form is 
the one on $L$ multiplied by $m$. 

Let $L$ be a non-degenerate lattice. Then the bilinear form of $L$ determines a 
canonical embedding $L\subset L^{\ast }=\Hom (L,\mathbb{Z})$. 
We denote by $A_{L}$ the factor group $L^{\ast }/L$ which is a finite abelian group.

Let $p$ be a prime number.
A lattice $L$ is called \textit{$p$-elementary} if $A_{L}\simeq (\mathbb{Z}/p\mathbb{Z})^{s}$,
where $s$ is the minimal number of generator of $A_{L}$.
For a $p$-elementary lattice we always have the inequality $s\leq r$, since 
$\mid L^{\ast }/L\mid =p^{s}$, $\mid L^{\ast }/pL^{\ast }\mid =p^{r}$ 
and $pL^{\ast }\subset L\subset L^{\ast }$. 

For example $A_{2}$, $E_{6}$, $E_{8}$, $E_{8}(3)$, $U$ and $U(3)$ are all $3$-elementary.

\vspace{3mm}

Even hyperbolic 3-elementary lattices were classified as follows.
\begin{thm}[\cite{R-S}]\label{3lattice}
An even, indefinite, $p$-elementary lattice of rank $n$ for $p \neq 2$ and $n\geq 2$ is 
uniquely determined by its discriminant (i.e., the number $s$).

For $p\neq 2$ a hyperbolic lattice corresponding to given value of $s\leq n$ exist if and only if
the following conditions are satisfied:
$n\equiv 0 \pmod 2$ and 

\begin{equation*}
\begin{cases}
\text{for } s\equiv 0 \pmod 2, & n\equiv 2 \pmod 4 \\
\text{for } s\equiv 1 \pmod 2, & p\equiv (-1)^{n/2-1} \pmod 4 
\end{cases}.
\end{equation*}

And moreover $n>s>0$, if $n\not \equiv 2 \pmod 8$.
\end{thm}

\vspace{3mm}

Let $L_{K3}$ be an even unimodular lattice of signature $(3,19)$. 
It is known that $L_{K3}$ is isometric to $U^{\oplus 3}\oplus E_{8}^{\oplus 2}$ 
by the classification of even unimodular indefinite lattices (\cite{Serre}, chapter 5, $\S 2$, Theorem 5).
Note that for a $K3$ surface $X$, $H^{2}(X, \mathbb{Z})$, together with the cup product, has a structure
of lattice which is isomorphic to $L_{K3}$.

The following proposition implies that 
if $X$ has non-symplectic automorphisms of order 3 which act trivially on $S_{X}$ 
then $22-\rho -2s\geq 0$.

\begin{prop}[\cite{Vorontsov}]\label{vo}
Assume that there exists a non-symplectic automorphism $\varphi $ of order $p$ on $X$ which acts trivially on $S_{X}$.
Then $S_{X}$ is a $p$-elementary lattice and $s\leq \rk T_{X}/(p-1)$.
\end{prop}
\begin{proof}
The same argument as in Lemma 1.3 in \cite{Oguiso-Zhang2} shows the assertion.
\end{proof}

Hence it is enough to classify even hyperbolic 3-elementary lattices 
admitting a primitive embedding in $L_{K3}$ satisfying $22-\rho -2s\geq 0$ 
(see \cite{Nikulin1} Sec.1, part $12^{\circ }$).

\begin{lem}\label{table}
Let $S$ be an even hyperbolic 3-elementary lattice admitting a primitive embedding in  $L_{K3}$. 
Let $T$ be the orthogonal complement of $S$ in $L_{K3}$. 
Then the next tables give all even hyperbolic 3-elementary lattices 
admitting a primitive embedding in $L_{K3}$ satisfying $22-\rk S -2s\geq 0$.

\begin{center}
\begin{longtable}{cc}

\begin{tabular}{|c|c|c|}
\hline
\multicolumn{3}{|c|}{$\rk S =2$} \\
\hline
\multicolumn{1}{|c|}{$s$} & \multicolumn{1}{c}{$S$} & \multicolumn{1}{|c|}{$T$} \\
\hline
0 & $U$ & $U^{\oplus 2}\oplus E_{8}^{\oplus 2}$ \\
\hline
2 & $U(3)$ & $U\oplus U(3)\oplus E_{8}^{\oplus 2}$  \\
\hline
\end{tabular}

& 

\begin{tabular}{|c|c|c|}
\hline
\multicolumn{3}{|c|}{$\rk S =4$} \\
\hline
\multicolumn{1}{|c|}{$s$}& \multicolumn{1}{c}{$S$} & \multicolumn{1}{|c|}{$T$} \\
\hline
1 & $U\oplus A_{2}$ & $U^{\oplus 2}\oplus E_{6}\oplus E_{8}$ \\
\hline
3 & $U(3)\oplus A_{2}$ & $U\oplus U(3)\oplus E_{6}\oplus E_{8}$  \\
\hline
\end{tabular}
\\
 
& \\
 
\begin{tabular}{|c|c|c|}
\hline
\multicolumn{3}{|c|}{$\rk S =6$} \\
\hline
\multicolumn{1}{|c|}{$s$}& \multicolumn{1}{c}{$S$} & \multicolumn{1}{|c|}{$T$} \\
\hline
2 & $U\oplus A_{2}^{\oplus 2}$ & $U^{\oplus 2}\oplus E_{6}^{\oplus 2}$ \\
\hline
4 & $U(3)\oplus A_{2}^{\oplus 2}$ & $U\oplus U(3)\oplus E_{6}^{\oplus 2}$  \\
\hline
\end{tabular}

& 

\begin{tabular}{|c|c|c|}
\hline
\multicolumn{3}{|c|}{$\rk S =8$} \\
\hline
\multicolumn{1}{|c|}{$s$}& \multicolumn{1}{c}{$S$} & \multicolumn{1}{|c|}{$T$} \\
\hline
1 & $U\oplus E_{6}$ & $U^{\oplus 2}\oplus E_{8}\oplus A_{2}$ \\
\hline
3 & $U\oplus A_{2}^{\oplus 3}$ & $U\oplus U(3)\oplus E_{8}\oplus A_{2}$  \\
\hline
5 & $U(3)\oplus A_{2}^{\oplus 3}$ & $A_{2}(-1)\oplus E_{6}\oplus A_{2}^{\oplus 3}$  \\
\hline
7 & $U(3)\oplus E_{6}^{\ast }(3)$ & $A_{2}(-1)\oplus A_{2}^{\oplus 6}$  \\
\hline
\end{tabular}
\\

& \\

\begin{tabular}{|c|c|c|}
\hline
\multicolumn{3}{|c|}{$\rk S =10$} \\
\hline
\multicolumn{1}{|c|}{$s$}& \multicolumn{1}{c}{$S$} & \multicolumn{1}{|c|}{$T$} \\
\hline
0 & $U\oplus E_{8}$ & $U^{\oplus 2}\oplus E_{8}$ \\
\hline
2 & $U\oplus E_{6}\oplus A_{2}$ & $U\oplus U(3)\oplus E_{8}$  \\
\hline
4 & $U\oplus A_{2}^{\oplus 4}$ & $U\oplus U(3)\oplus E_{6}\oplus A_{2}$  \\
\hline
6 & $U(3)\oplus A_{2}^{\oplus 4}$ & $A_{2}(-1)\oplus A_{2}^{\oplus 5}$  \\
\hline
\end{tabular}

& 

\begin{tabular}{|c|c|c|}
\hline
\multicolumn{3}{|c|}{$\rk S =12$} \\
\hline
\multicolumn{1}{|c|}{$s$}& \multicolumn{1}{c}{$S$} & \multicolumn{1}{|c|}{$T$} \\
\hline
1 & $U\oplus E_{8}\oplus A_{2}$ & $A_{2}(-1)\oplus E_{8}$ \\
\hline
3 & $U\oplus E_{6}\oplus A_{2}^{\oplus 2}$ & $A_{2}(-1)\oplus E_{6}\oplus A_{2}$  \\
\hline
5 & $U\oplus A_{2}^{\oplus 5}$ & $A_{2}(-1)\oplus A_{2}^{\oplus 4}$  \\
\hline
\end{tabular}
\\

& \\

\begin{tabular}{|c|c|c|}
\hline
\multicolumn{3}{|c|}{$\rk S =14$} \\
\hline
\multicolumn{1}{|c|}{$s$}& \multicolumn{1}{c}{$S$} & \multicolumn{1}{|c|}{$T$} \\
\hline
2 & $U\oplus E_{8}\oplus A_{2}^{\oplus 2}$ & $A_{2}(-1)\oplus E_{6}$ \\
\hline
4 & $U\oplus E_{6}\oplus A_{2}^{\oplus 3}$ & $A_{2}(-1)\oplus A_{2}^{\oplus 3}$  \\
\hline
\end{tabular}

&

\begin{tabular}{|c|c|c|}
\hline
\multicolumn{3}{|c|}{$\rk S =16$} \\
\hline
\multicolumn{1}{|c|}{$s$}& \multicolumn{1}{c}{$S$} & \multicolumn{1}{|c|}{$T$} \\
\hline
1 & $U\oplus E_{8}\oplus E_{6}$ & $U^{\oplus 2}\oplus A_{2}$ \\
\hline
3 & $U\oplus E_{8}\oplus A_{2}^{\oplus 3}$ & $A_{2}(-1)\oplus A_{2}^{\oplus 2}$  \\
\hline
\end{tabular}
\\

& \\

\begin{tabular}{|c|c|c|}
\hline
\multicolumn{3}{|c|}{$\rk S =18$} \\
\hline
\multicolumn{1}{|c|}{$s$}& \multicolumn{1}{c}{$S$} & \multicolumn{1}{|c|}{$T$} \\
\hline
0 & $U\oplus E_{8}^{\oplus 2}$ & $U^{\oplus 2}$ \\
\hline
2 & $U\oplus E_{8}\oplus E_{6}\oplus A_{2}$ & $U\oplus U(3)$  \\
\hline
\end{tabular}

& 

\begin{tabular}{|c|c|c|}
\hline
\multicolumn{3}{|c|}{$\rk S =20$} \\
\hline
\multicolumn{1}{|c|}{$s$}& \multicolumn{1}{c}{$S$} & \multicolumn{1}{|c|}{$T$} \\
\hline
1 & $U\oplus E_{8}^{\oplus 2}\oplus A_{2}$ & $A_{2}(-1)$ \\
\hline
\end{tabular}
\\
\caption[]{}\label{3el}
\end{longtable}
\end{center}
\end{lem}

\begin{proof}
Since $L_{K3}$ is a unimodular lattice and $S$ is a primitive sublattice of $L_{K3}$, $A_{S}\simeq A_{T}$ (\cite{BHPV} Lemma 2.5). 
By Theorem \ref{3lattice}, 
if $\rk S=2, 6, 10, 14, 18$ then $s \equiv 0 \pmod 2$ and 
if $\rk S=4, 8, 12, 16, 20$ then $s \equiv 1 \pmod 2$. 
Moreover if $\rk S=6, 14$ then $s > 0$.
\end{proof}

\section{The fixed locus of non-symplectic automorphisms}\label{form}

In this section, we shall see that the fixed locus 
$X^{\varphi }$ is determined by the invariants of the N\'{e}ron-Severi lattice $S_{X}$.

\begin{lem}\label{sayou}
Let $X$ be an algebraic $K3$ surface, $\varphi $ a non-symplectic automorphism of order $3$ on $X$.
Then we have :
\begin{itemize}
\item[(1)] $\varphi ^{\ast }\mid T_{X}\otimes \mathbb{C}$ can be diagonalized as:
\[ \begin{pmatrix} \zeta I_{r} & 0 \\ 0 & \bar{\zeta }I_{r}  \end{pmatrix} \]
where $I_{r}$ is the identity matrix of size $r$, $\zeta$ is a primitive third root of unity.
\item[(2)] Let $P$ be an isolated fixed point of $\varphi $ on $X$. Then 
$\varphi ^{\ast }$ can be written as 
\[ \begin{pmatrix}  \zeta ^{2} & 0 \\ 0 & \zeta ^{2}  \end{pmatrix} \]
under some appropriate local coordinates around $P$.
\item[(3)] Let $C$ be a fixed irreducible curve and $Q$ a point on $C$. 
Then $\varphi ^{\ast }$ can be written as
\[ \begin{pmatrix}  1 & 0 \\ 0 & \zeta   \end{pmatrix} \] 
under some appropriate local coordinates around $Q$. 
In particular, fixed curves are non-singular.
\end{itemize}
\end{lem}

\begin{proof}
(1) This follows form \cite{Nikulin2}, Theorem 3.1.  

(2), (3) 
Since $\varphi ^{\ast }$ acts on $H^{0}(X, \Omega _{X}^{2})$ as a multiplication by $\zeta $, 
it acts on the tangent space of a fixed point as 
\[ \begin{pmatrix} 1 & 0 \\ 0 & \zeta  \end{pmatrix} \hspace{1cm} \text{or}  \hspace{1cm} \begin{pmatrix} \zeta^{2} & 0 \\ 0 & \zeta^{2}  \end{pmatrix}. \]
\end{proof}
Thus the fixed locus of $\varphi $ consists of disjoint union of non-singular curves and isolated points. 
Hence we can express the irreducible decomposition of $X^{\varphi }$ as 
\[ X^{\varphi } =\{ P_{1} \} \amalg \dots \amalg \{ P_{M} \} \amalg C_{1} \amalg \dots \amalg C_{N}, \] 
where $P_{j}$ is an isolated point and $C_{k}$ is a non-singular curve.

\begin{prop}\label{point}
Let $\rho $ be the Picard number of $X$. Then the number of isolated points $M$ is $\rho /2-1$.
\end{prop}
\begin{proof}
First we calculate the holomorphic Lefschetz number $L(\varphi )$ in two ways as in 
\cite{AS1}, page 542 and \cite{AS2}, page 567. That is 
\[ L(\varphi ) = \sum _{i=0}^{2} \text{tr}(\varphi ^{\ast }|H^{i}(X, \mathcal{O}_{X})), \]
\[ L(\varphi ) = \sum _{j=1}^{M} a(P_{j}) + \sum _{k=1}^{N}b(C_{k}). \]

Here 
\begin{align*}
a(P_{j}) : & =\frac{1}{\det (1-\varphi ^{\ast }|T_{P_{j}})}  \\
 & =\frac{1}{\det \left ( \begin{pmatrix} 1 & 0 \\ 0 & 1 \end{pmatrix} -  \begin{pmatrix} \zeta ^{2} & 0 \\ 0 & \zeta ^{2} \end{pmatrix} \right ) } \\
 & =-\frac{\zeta}{3}, 
\end{align*}
\begin{align*}
b(C_{k}) : & =\frac{1-g(C_{k})}{1-\zeta ^{-2}}- \frac{\zeta ^{-2}C_{k}^{2}}{(1-\zeta ^{-2})^{2}} \\
 & =\frac{\zeta (1-g(C_{k}))}{3}, 
\end{align*}
where $T_{P_{j}}$ is the tangent space of $X$ at $P_{j}$, $g(C_{k})$ is the genus of $C_{k}$ and 
$\zeta ^{2}$ is the eigenvalue of the action $\varphi _{\ast }$ on the normal bundle of $C_{k}$ (cf. Lemma \ref{sayou}).

Using the Serre duality $H^{2}(X, \mathcal{O}_{X})\simeq H^{0}(X,\mathcal{O}_{X}(K_{X}))^{\vee }$, we calculate 
from the first formula that $L(\varphi )=1+\zeta ^{2}=-\zeta $. From the second formula, we obtain
\[ L(\varphi ) =\frac{-M\zeta }{3} + \sum_{k=1}^{N} \frac{\zeta (1-g(C_{k})) }{3}. \]

Combing these two formulae, we have 
\begin{equation}\label{3}
M-\sum _{k=1}^{N}(1-g(C_{k}))=3.
\end{equation}
Next we apply the topological Lefschetz formula:
\[ \chi _{\text{top}}(X^{\varphi }) =\sum _{i=0}^{4}(-1)^{i} \text{tr}(\varphi ^{\ast }|H^{i}(X, \mathbb{R})). \]
The left-hand side is 
\begin{equation}\label{left}
\chi _{\text{top}}(X^{\varphi })=M+\sum _{k=1}^{N}(2-2g(C_{k})). 
\end{equation}

By Lemma \ref{sayou}, one has the following diagonalized actions
\[ \varphi ^{\ast }|T_{X}\otimes \mathbb{C}= \begin{pmatrix} \zeta I_{r} & 0 \\ 0 & \bar{\zeta }I_{r}  \end{pmatrix}. \]

Since $\varphi ^{\ast }$ acts trivially on $S_{X}$, $\text{tr}(\varphi ^{\ast }|S_{X})=\rho $.
Thus we have $\rho +2r=22$. 
Hence we can calculate the right -hand side of the Lefschetz formula as follows:

\begin{align}\label{right}
\sum _{i=0}^{4}(-1)^{i} \text{tr}(\varphi ^{\ast }|H^{i}(X, \mathbb{R})) : & = 1-0+\text{tr}(\varphi ^{\ast }|S_{X})+\text{tr}(\varphi ^{\ast }|T_{X})-0+1 \notag  \\
& = 2+\rho -r \notag \\
& = 2+\rho -\frac{22-\rho }{2}  \notag \\
& = \frac{3\rho -18}{2}.
\end{align}

By (\ref{3}), (\ref{left}) and (\ref{right}), $M=\rho /2-1$.
\end{proof}

By the Hodge index theorem, the following three cases are possible:
\begin{equation}\label{nashi}
X^{\varphi } = \phi \tag{A}; 
\end{equation}
\begin{equation}\label{genus1}
X^{\varphi } = 
\{ P_{1} \} \amalg \dots \amalg \{ P_{M} \} 
\amalg C_{1}^{(1)} \amalg \dots \amalg C_{L}^{(1)}
\amalg E_{1} \amalg \dots \amalg E_{K} \tag{B}; 
\end{equation}
\begin{equation}\label{general}
X^{\varphi } = \{ P_{1} \} \amalg \dots \amalg \{ P_{M} \} \amalg C^{(g)} \amalg E_{1} \amalg \dots \amalg E_{N-1} \tag{C}, 
\end{equation}
where we denote by $P_{i}$ an isolated point, $C^{(g)}$ a non-singular curve of genus $g$ and by $E_{j}$ a non-singular rational curve.

\vspace{3mm}

The following lemma follows from \cite{PSS} $\S 3$ Corollary 3 and the 
classification of singular fibers of elliptic fibrations \cite{Kodaira}. 

\begin{lem}\label{elliptic}
Let $X$ be an algebraic $K3$ surface.
Assume that $S_{X}=U(m)\oplus K_{1}\oplus \dots \oplus K_{r}$, 
where $m=$ 1 or 3, and $K_{i}$ is a lattice isomorphic to $A_{2}$, $E_{6}$ or $E_{8}$.
Then there exists an elliptic fibration $\pi :X\longrightarrow \mathbb{P}^{1}$.
Moreover $\pi $ has a reducible singular fiber whose dual graph is of type $\widetilde{K_{i}}$.
\end{lem}

\begin{rem}\label{elpss}
Let $ \{ e, f \}$ be a basis of $U$ (resp. $U(3)$) with $\langle e, e \rangle = \langle f, f \rangle =0$ 
and $\langle e, f \rangle =1$ (resp. $\langle e, f \rangle =3$ ) .
If necessary replacing $e$ by $\varphi (e)$, 
where $\varphi $ is a composition of reflections induced from non-singular 
rational curves on $X$, we may assume that $e$ is represented by the class of an elliptic curve $F$ 
and the linear system $|F|$ defines an elliptic fibration $\pi :X\longrightarrow \mathbb{P}^{1}$.
Note that $\pi $ has a section $e-f$ in case $U$. 
In case $U(3)$, there are no $(-2)$-vectors $r$ with $\langle r, e \rangle=1$, and hence $\pi$ has no sections.
\end{rem}

It follows from Remark \ref{elpss} and Table \ref{3el} that 
$X$ has an elliptic fibration $\pi :X\longrightarrow \mathbb{P}^{1}$.
In the following, we fix such the elliptic fibration.

\begin{lem}\label{246}
In the case (\ref{general}), if  $\rho < 8$ then the type of singular fibers is II or IV.
\end{lem}
\begin{proof}
We obtain
$\chi _{\text{top}} \left (\sum P_{j} \right ) =\rho/2-1$ and $\chi _{\text{top}}(X^{\varphi })=(3\rho -18)/2$
by Proposition \ref{point} and the topological Lefschetz formula.

From these two equations, we get 
\begin{equation}\label{euler}
\chi _{\text{top}} \left (C^{(g)} \amalg \sum_{k=1}^{N-1} E_{k} \right ) =\rho -8. 
\end{equation}
The left-hand side is calculated as  
\begin{equation}\label{euler-left}
\chi _{\text{top}} \left (C^{(g)} \amalg \sum_{k=1}^{N-1} E_{k} \right ) =(2-2g)+2(N-1).
\end{equation}
Therefore if $\rho < 8$ then there exists a non-singular curve $C^{(g)}$ with $g\geq 2$.

Since $S_{X}=U(m)\oplus A_{2}^{\oplus l}$, where $m=$1 or 3, $0\leq l\leq 3$, 
the type of singular fibers of $\pi $ is I$_{1}$, II, I$_{3}$ or IV.
Let $F$ be a singular fiber of $\pi $.  
By the Hodge index theorem, the intersection number $C^{(g)}.F$ is positive. 
This implies that the automorphism $\varphi $ acts trivially on the base.
Hence a smooth fiber has an automorphism of order 3.
Thus a fixed locus of a smooth fiber is exactly three isolated points.
Therefore the functional invariant of $\pi $ is 0.
This implies that if  $\rho < 8$ then the type of singulars fiber is II or IV.
\end{proof}

First we take up the case (\ref{nashi}) and the case (\ref{genus1}).

\begin{lem}
The case (\ref{nashi}) does not occur, that is, $X^{\varphi } \neq \phi$.
\end{lem}
\begin{proof}
By Proposition \ref{point}, if $X^{\varphi }=\phi  $ then  $\rho =2$.  
But if $\rho =2$ then $\chi_{\text{top}} (X^{\varphi })=-6$ by the topological Lefschetz formula. 
This is a contradiction by $\chi_{\text{top}} (\phi )=0$.
Hence the case (\ref{nashi}) does not occur.
\end{proof}

\begin{lem}\label{c23}
In the case (\ref{genus1}), $L\leq 1$. 
\end{lem}
\begin{proof}
An elliptic curve $C_{j}^{(1)}$ belongs to one elliptic pencil $|C_{j}^{(1)}|:X\longrightarrow \mathbb{P}^{1}$.
Assume $L\geq 2$ and $K>0$. 
Then $\varphi $ fixes at least three fibers: $C_{1}^{(1)},\dots ,C_{L}^{(1)}$ and a fiber containing $E_{1}$.
Since an automorphism of order 3 on $\mathbb{P}^{1}$ has exactly two isolated fixed points, 
$\varphi $ is trivial on the base $\mathbb{P}^{1}$. 
And it is also trivial on a fiber $C_{j}^{(1)}$. 
Hence $\varphi $ is a symplectic automorphism. 
This is a contradiction. 

We remark that if $L=2$ and $K=0$ then $M>0$ by Proposition \ref{point}.
Actually if $\rho =2$ then $\chi _{\text{top}}(X^{\varphi })=-6$ by the topological Lefschetz formula.
This implies that $X^{\varphi }$ contains a non-singular curve $C^{(g)}$ with genus $g\geq 2$.
This is a contradiction by the Hodge index theorem.
Hence $\varphi $ fixes at least three fibers: $C_{1}^{(1)}, C_{2}^{(1)}$ and a fiber containing $P_{1}, \dots ,P_{M}$.
This implies that $\varphi $ is trivial on the base $\mathbb{P}^{1}$.
Therefore the case of $L=2$ and $K=0$ does not occur.
\end{proof}

The Lemma \ref{c23} implies that the case (\ref{genus1}) is a special case of the case (\ref{general}).
Hence more generally, we take up the case (\ref{general}). 
Actually we have the following results.

\begin{thm}\label{mt1}
Let $S_{X}$ be the N\'{e}ron-Severi lattice except $U(3)\oplus E_{6}^{\ast }(3)$, 
let $\rho$ be the Picard number of $X$ and let $s$ be the minimal number of generators of $S_{X}^{\ast }/S_{X}$.

If $22-\rho -2s\geq 0$, then 
$X$ has a non-symplectic automorphism $\varphi $ of order 3 which acts trivially on $S_{X}$. 
Moreover the fixed locus $X^{\varphi }$ has the form 
\[ X^{\varphi }=\{ P_{1} \} \amalg \dots \amalg \{ P_{M} \} \amalg C^{(g)} \amalg E_{1} \amalg \dots \amalg E_{N-1} \]
where $M = \rho /2-1$, $g=(22-\rho -2s)/4$, $N=(6+\rho -2s)/4$.
\end{thm}

\begin{proof}

The equation $M = \rho /2-1$ follows from Proposition \ref{point}.
We shall prove $g=(22-\rho -2s)/4$ and $N=(6+\rho -2s)/4$.

\vspace{3mm}

\noindent (Case 1) $\rho < 8$.

Let $F$ be a fiber of $\pi $. 
From the proof of Lemma \ref{246}, $\varphi $ acts trivially on the base of $\pi $.
Thus if $\pi $ has a section (resp. no section) then $C^{(g)}.F=2$ (resp. $C^{(g)}.F=3$) in this case.

It is known that 
\begin{equation}\label{fiber24}
\sum _{\text{$F$: singular fiber}} \chi _{\text{top}}(F) = \chi _{\text{top}}(X) =24. 
\end{equation}
Since the Euler number of a singular fiber of type II or of type IV is 2 or 4, respectively, 
the number of singular fibers of type II is $[24-4\{(\rho -2)/2\}]/2=14-\rho $ by equation (\ref{fiber24}).
If $\pi $ has a section (resp. no section) then the number of singular fibers of type IV $=s$ (resp. $s-2$).
Thus $(\rho -2)/2=s$ (resp. $(\rho -2)/2=s-2$).

We remark that $C^{(g)}$ meets a singular fiber of type II at cusp.
Since if $\pi $ has a section then $C^{(g)}.F=2$, by the Hurwitz formula, we have 
\begin{align*}
2g-2 & =2(2g(\mathbb{P}^{1})-2)+(14-\rho )(2-1) \\
2g    & =\frac{24-2\rho }{2}\\ 
       & =\frac{24-\rho -(2s-2)}{2}\\
g      & =\frac{22-\rho -2s}{4}.
\end{align*}

Next we assume $\pi $ has no section.

By Lemma \ref{246}, if $X^{\varphi }$ contains some non-singular rational curves $E_{j}$ then 
these are components of singular fibers of type IV.
But this is a contradiction by $C^{(g)}.F=3$. 
Hence $X^{\varphi }$ contains no non-singular rational curves. 
By the Lefschetz formula, 
\begin{align*}
(2-2g) + \frac{\rho -2}{2} & =\frac{3\rho -18}{2} \\
2g                                & =\frac{20-2\rho }{2}\\ 
                                   & =\frac{20-\rho -(2s-2)}{2}\\
g                                 & =\frac{22-\rho -2s}{4}.
\end{align*}

Moreover by this equation, (\ref{euler}) and (\ref{euler-left}), we have 
\[ N=\frac{6+\rho -2s}{4}. \]

\vspace{3mm}

\noindent (Case 2) $\rho \geq 8$ and $\pi $ has a section.

We will check the formulae $g=(22-\rho -2s)/4$ and $N=(6+\rho -2s)/4$ for each $S_{X}$ individually.

Assume $S_{X}=U\oplus E_{8}\oplus E_{6}$. 
By (\ref{euler}) and (\ref{euler-left}), $g=N-4$. 
Since $\rk S_{X}=16$, we have $\chi _{\text{top}} \left (C^{(g)} \amalg \sum_{k=1}^{N-1} E_{k} \right ) =8$.
Thus the number of non-singular rational curves contained in $X^{\varphi }$ is 4 or more.

Note $\pi $ has singular fibers of type II$^{\ast }$ and of type IV$^{\ast }$. 
If $g=0$ (i.e. $N=4$) then the number of isolated fixed points of $\varphi $ is exactly eight. 
This is because the component with multiplicity 6 of singular fiber of type II$^{\ast }$ 
and the component with multiplicity 3 of singular fiber of type IV$^{\ast }$
are pointwisely fixed by $\varphi $.
This is a contradiction by Proposition \ref{point}.

If $g\geq 2$ then the automorphism $\varphi $ acts trivially on the base of an elliptic fibration. 
Since the Euler number of  a singular fiber of type II$^{\ast }$ and of type IV$^{\ast }$ is 10 or 8, respectively, 
$\pi $ has three singular fibers of type II. 
By the Hurwitz formula, we have 
\[ 2g-2=2(2g(\mathbb{P}^{1})-2)+4(2-1). \]
But this is a contradiction by $g\geq 2$.
Hence $g=1$ and $N=5$. 
Therefore $X^{\varphi }:=C^{(1)}\amalg \coprod_{i=1}^{4} \mathbb{P}^{1}_{i}\amalg \coprod_{j=1}^{7} \{P_{j}\}$.

Similarly in other cases we can calculate the genus and the number of the fixed curves 
by the same argument of these examples.  
These results satisfy the assertion.

\vspace{3mm}

\noindent (Case 3) $\rho \geq 8$ and $\pi $ has no sections.

Assume $S_{X}=U(3)\oplus A_{2}^{\oplus 4}$. By (\ref{euler}) and (\ref{euler-left}), $g=N-1$. 
Since $\rk S_{X}=10$, $\chi _{\text{top}} \left (C^{(g)} \amalg \sum_{k=1}^{N-1} E_{k} \right ) =2$.
Thus the number of non-singular rational curves which $X^{\varphi }$ contains is one or more.

Note $\pi $ has four singular fibers $F_{i}$ of type IV or of type I$_{3}$. 
Thus the automorphism $\varphi $ acts trivially on the base of an elliptic fibration.
Since $\pi $ has no section, $C^{(g)}.F_{i}=3$.
Moreover we can see that all $F_{i}$ are of type IV by Lemma \ref{elliptic}.
And a singular fiber of type IV has exactly one isolated fixed point at the intersection point of 
each component of the singular fiber. 
Therefore $\varphi $ has exactly one fixed curve. 
Hence $X^{\varphi }:=C^{(0)}\amalg \{P_{1}\} \amalg \{P_{2}\} \amalg \{P_{3}\}\amalg \{P_{4}\}$.

A similar assertion holds in the other cases.
\end{proof}

Finally we consider the case of $S_{X}=U(3)\oplus E_{6}^{\ast }(3)$.

\begin{prop}\label{mt2}
If $S_{X}=U(3)\oplus E_{6}^{\ast }(3)$ then 
$X^{\varphi }=\{ P_{1} \} \amalg \{ P_{2} \} \amalg \{ P_{3} \}$.
\end{prop}
\begin{proof}
First we show that $S_{X}$ has no $(-2)$ vectors. 
Let $q_{E_{6}}$ be a discriminant form of $E_{6}$ and let $\bar{a}$ be a generator of  $E_{6}^{\ast }/E_{6}$. 
Since $U\oplus E_{6}^{\ast }/U\oplus E_{6}\simeq E_{6}^{\ast }/E_{6}$ and $q_{E_{6}}(\bar{a})=-4/3$, 
$U\oplus E_{6}^{\ast }$ has no $(-2/3)$ vectors. 
This implies that $S_{X}$ has no $(-2)$ vectors.
Hence $X$ has no non-singular rational curves, i.e. $K=0$.

Since $\rho =8$, we have $M=3$ by Proposition \ref{point}. 
From (\ref{euler}), (\ref{euler-left}) and Lemma \ref{c23}, if there exits a fixed curve $C^{(g)}$ then $g=1$.

Now we assume that $X^{\varphi }$ contains a fixed curve, 
i.e. $X^{\varphi }=\{P_{1}\}\amalg \{P_{2}\} \amalg \{P_{3}\} \amalg C^{(1)}$.
Let $\{e, f\}$ be a basis of $U(3)$. 
By Remark \ref{elpss}, $|e|$ define an elliptic fibration $\pi:X\longrightarrow \mathbb{P}^{1}$. 
We assume first $\varphi $ acts trivially on the base of $\pi $. 
Then every smooth fiber of $\pi $ has an automorphisms of order 3, 
and hence the functional invariant is 0. 
Moreover $C^{(1)}$ meets a smooth fiber at three points. 
If $\pi $ has a reducible fiber then $X$ has a non-singular rational curve as a component of the reducible fiber. 
This is a contradiction. 
Hence $\pi $ has 12 singular fibers of type II by (\ref{fiber24}). 
Obviously 12 cusps of 12 singular fibers are contained in $X^{\varphi}$. 
Hence at least 9 cusps lie on $C^{(1)}$. 
Then $\pi|_{C^{(1)}}:C^{(1)}\longrightarrow \mathbb{P}^{1}$ is a covering of degree 3 
ramified at these 9 cusps. 
On the other hand, the Hurwitz formula implies that 
\[0=2g(C^{(1)})-2\geq 3(2g(\mathbb{P}^{1})-2)+9(2-1)=3.\] 
This is a contradiction.

Next we assume that $\varphi $ acts on the base of $\pi $ as an automorphism of order $3$.
Then $\varphi $ has exactly two isolated fixed points $Q_{1}$ and $Q_{2}$ on the base of $\pi $.
Then $C^{(1)}$ is equal to $\varphi ^{-1}(Q_{1})$ or $\varphi ^{-1}(Q_{2})$.
Since $C^{(1)}.F=e.f=3$ and $C^{(1)}$ is fixed by $\varphi $, 
$\varphi $ acts trivially on the base of $|f|$ where $F$ is a smooth fiber of $|f|$.
Hence this case reduces to the argument as above. 

Thus $X^{\varphi }$ has no fixed curve.
Hence $X^{\varphi }=\{ P_{1} \} \amalg \{ P_{2} \} \amalg \{ P_{3} \}$.
\end{proof}

Theorem \ref{mt} follows from Proposition \ref{vo}, Theorem \ref{mt1} and Proposition \ref{mt2}.

\begin{rem}
The case of $S_{X}=U(3)\oplus E_{6}^{\ast }(3)$ satisfies the equation $N=(6+\rho -2s)/4$.
\end{rem}

Finally we give the list of fixed loci in Table \ref{SandF}.

\begin{longtable}{|c|c|}
\hline
$S_{X}$ & $X^{\varphi }$ \\
\hline
$U$ & $C^{(5)}\amalg \mathbb{P}^{1}$ \\
\hline
$U(3)$ & $C^{(4)}$ \\
\hline
$U\oplus A_{2}$ & $C^{(4)}\amalg \mathbb{P}^{1}\amalg \{pt\}$ \\
\hline
$U(3)\oplus A_{2}$ & $C^{(3)}\amalg \{pt\}$ \\
\hline
$U\oplus A_{2}^{\oplus 2}$ & $C^{(3)}\amalg \mathbb{P}^{1}\amalg \{pt\} \times 2$ \\
\hline
$U(3)\oplus A_{2}^{\oplus 2}$ & $C^{(2)}\amalg \{pt\} \times 2$ \\
\hline
$U\oplus E_{6}$ & $C^{(3)}\amalg \mathbb{P}^{1}\times 2 \amalg \{pt\} \times 3$ \\
\hline
$U\oplus A_{2}^{\oplus 3}$ & $C^{(2)}\amalg \mathbb{P}^{1}\amalg \{pt\} \times 3$ \\
\hline
$U(3)\oplus A_{2}^{\oplus 3}$ & $C^{(1)}\amalg \{pt\} \times 3$ \\
\hline
$U(3)\oplus E_{6}^{\ast }(3)$ & $\{pt\} \times 3$ \\
\hline
$U\oplus E_{8}$ & $C^{(3)}\amalg \mathbb{P}^{1}\times 3\amalg \{pt\} \times 4$ \\
\hline
$U\oplus E_{6}\oplus A_{2}$ & $C^{(2)}\amalg \mathbb{P}^{1}\times 2\amalg \{pt\} \times 4$ \\
\hline
$U\oplus A_{2}^{\oplus 4}$ & $C^{(1)}\amalg \mathbb{P}^{1}\amalg \{pt\} \times 4$ \\
\hline
$U(3)\oplus A_{2}^{\oplus 4}$ & $C^{(0)}\amalg \{pt\} \times 4$ \\
\hline
$U\oplus E_{8}\oplus A_{2}$ & $C^{(2)}\amalg \mathbb{P}^{1}\times 3\amalg \{pt\} \times 5$ \\
\hline
$U\oplus E_{6}\oplus A_{2}^{\oplus 2}$ & $C^{(1)}\amalg \mathbb{P}^{1}\times 2\amalg \{pt\} \times 5$ \\
\hline
$U\oplus A_{2}^{\oplus 5}$ & $C^{(0)}\amalg \mathbb{P}^{1}\amalg \{pt\} \times 5$ \\
\hline
$U\oplus E_{8}\oplus A_{2}^{\oplus 2}$ & $C^{(1)}\amalg \mathbb{P}^{1}\times 3\amalg \{pt\} \times 6$ \\
\hline
$U\oplus E_{6}\oplus A_{2}^{\oplus 3}$ & $C^{(0)}\amalg \mathbb{P}^{1}\times 2\amalg \{pt\} \times 6$ \\
\hline
$U\oplus E_{8}\oplus E_{6}$ & $C^{(1)}\amalg \mathbb{P}^{1}\times 4\amalg \{pt\} \times 7$ \\
\hline
$U\oplus E_{8}\oplus A_{2}^{\oplus 3}$ & $C^{(0)}\amalg \mathbb{P}^{1}\times 3\amalg \{pt\} \times 7$ \\
\hline
$U\oplus E_{8}^{\oplus 2}$ & $C^{(1)}\amalg \mathbb{P}^{1}\times 5\amalg \{pt\} \times 8$ \\
\hline
$U\oplus E_{8}\oplus E_{6}\oplus A_{2}$ & $C^{(0)}\amalg \mathbb{P}^{1}\times 4\amalg \{pt\} \times 8$ \\
\hline
$U\oplus E_{8}^{\oplus 2}\oplus A_{2}$ & $C^{(0)}\amalg \mathbb{P}^{1}\times 5\amalg \{pt\} \times 9$ \\
\hline
\caption[]{N\'{e}ron-Severi lattices and fixed loci}\label{SandF}
\end{longtable}

\section{Examples}\label{exa}

In this section, we give examples of algebraic $K3$ surfaces with a non-symplectic automorphism of order 3. 
We remark that a $K3$ surface $X$ has an elliptic fibration $\pi :X\longrightarrow \mathbb{P}^{1}$ 
from Remark \ref{elpss} and Table \ref{3el}.
First we consider the case that $\pi $ has a section.

\begin{example}\label{aek3}
In the following we give affine equations of elliptic $K3$ surfaces. 
We define an automorphism $\varphi $ of $X$ as follows: $\varphi (x, y, z, u )=(x, y, \zeta z, u )$ 
where $\zeta $ is a primitive third root of unity.

\begin{longtable}{|c|c|}
\hline
$S_{X}$ & definition equation \\
\hline
$U$ & $\displaystyle  z^{3}=y \left( y^{2}\prod_{i=1}^{12}(u-a_{i})-x^{2} \right) $ \\
\hline
$U\oplus A_{2}$ & $\displaystyle  z^{3}=y \left( y^{2}(u-a_{0})^{2}\prod_{i=1}^{10}(u-a_{i})-x^{2} \right) $ \\
\hline
$U\oplus A_{2}^{\oplus 2}$ & $\displaystyle  z^{3}=y \left( y^{2}(u-a_{0})^{2}\prod_{i=1}^{8}(u-a_{i})(u-a_{9})^{2}-x^{2} \right) $ \\
\hline
$U\oplus E_{6}$ & $\displaystyle  z^{3}=y \left( y^{2}(u-a_{0})^{4}\prod_{i=1}^{8}(u-a_{i})-x^{2} \right) $ \\
\hline
$U\oplus A_{2}^{\oplus 3}$ & $\displaystyle  z^{3}=y \left( y^{2}\prod_{i=1}^{6}(u-a_{i}) \prod_{j=7}^{9}(u-a_{j})^{2}-x^{2} \right) $ \\
\hline
$U\oplus E_{8}$ & $\displaystyle  z^{3}=y \left( y^{2}(u-a_{0})^{5}\prod_{i=1}^{7}(u-a_{i})-x^{2} \right) $ \\
\hline
$U\oplus E_{6}\oplus A_{2}$ & $\displaystyle  z^{3}=y \left( y^{2}(u-a_{0})^{4}\prod_{i=1}^{6}(u-a_{i})(u-a_{7})^{2}-x^{2} \right) $ \\
\hline
$U\oplus A_{2}^{\oplus 4}$ & $\displaystyle  z^{3}=y \left( y^{2}\prod_{i=1}^{4}(u-a_{i})\prod_{j=5}^{8}(u-a_{j})^{2}-x^{2} \right) $ \\
\hline
$U\oplus E_{8}\oplus A_{2}$ & $\displaystyle  z^{3}=y \left( y^{2}(u-a_{0})^{5}\prod_{i=1}^{5}(u-a_{i})(u-a_{6})^{2}-x^{2} \right) $ \\
\hline
$U\oplus E_{6}\oplus A_{2}^{\oplus 2}$ & $\displaystyle  z^{3}=y \left( y^{2}(u-a_{0})^{4}\prod_{i=1}^{4}(u-a_{i})\prod_{i=5}^{6}(u-a_{i})^{2}-x^{2} \right) $ \\
\hline
$U\oplus A_{2}^{\oplus 5}$ & $\displaystyle  z^{3}=y \left( y^{2}\prod_{i=1}^{2}(u-a_{i})\prod_{i=3}^{7}(u-a_{i})^{2}-x^{2} \right) $ \\
\hline
$U\oplus E_{8}\oplus A_{2}^{\oplus 2}$ & $\displaystyle  z^{3}=y \left( y^{2}(u-a_{0})^{5}\prod_{i=1}^{3}(u-a_{i})\prod_{i=4}^{5}(u-a_{i})^{2}-x^{2} \right) $ \\
\hline
$U\oplus E_{6}\oplus A_{2}^{\oplus 3}$ & $\displaystyle  z^{3}=y \left( y^{2}(u-a_{0})^{4}\prod_{i=1}^{3}(u-a_{i})^{2}\prod_{j=4}^{5}(u-a_{j})-x^{2} \right) $ \\
\hline
$U\oplus E_{8}\oplus E_{6}$ & $\displaystyle z^{3}=y \left( y^{2}(u-a_{0})^{5}\prod_{i=1}^{3}(u-a_{i})(u-a_{4})^{4}-x^{2} \right) $ \\
\hline
$U\oplus E_{8}\oplus A_{2}^{\oplus 3}$ & $\displaystyle z^{3}=y \left( y^{2}(u-a_{0})^{5}\prod_{i=1}^{3}(u-a_{i})^{2}(u-a_{4})-x^{2} \right) $ \\
\hline
$U\oplus E_{8}^{\oplus 2}$ & $\displaystyle z^{3}=y \left( y^{2}\prod_{i=1}^{2}(u-a_{i})^{5}\prod_{j=3}^{4}(u-a_{j})-x^{2} \right) $ \\
\hline
$U\oplus E_{8}\oplus E_{6}\oplus A_{2}$ & $z^{3}=y \left( y^{2}(u-a_{0})^{5}(u-a_{1})^{4}(u-a_{2})^{2}(u-a_{3})-x^{2} \right) $ \\
\hline
$U\oplus E_{8}^{\oplus 2}\oplus A_{2}$ & $\displaystyle z^{3}=y \left( y^{2}(u-a_{0})^{2}\prod_{i=1}^{2}(u-a_{i})^{5}-x^{2} \right) $ \\
\hline
\caption[]{N\'{e}ron-Severi lattices and definition equations}
\end{longtable}
\end{example}

Next we study one case in Example \ref{aek3} in detail. 

\textbf{Case: $S_{X}=U$} (\cite{Kondo1} (2.1)) \ 
Let $[x:y:z]$ be a system of a homogeneous coordinate of $\mathbb{P}^{2}$.
We take two copies $W_{0}:=\mathbb{P}^{2}\times \mathbb{C}_{0}$ and $W_{1}:=\mathbb{P}^{2}\times \mathbb{C}_{1}$ 
of the cartesian product $\mathbb{P}^{2}\times \mathbb{C}$ and form their union $W=W_{0}\cup W_{1}$ 
by identifying$([ x:y:z], u) \in W_{0}$ with $([ x_{1}:y_{1}:z_{1}], u_{1}) \in W_{1}$ if and only if 
$u=1/u_{1}$, $x=x_{1}$, $u=u_{1}^{6}y_{1}$ and $z=u_{1}^{2}z_{1}$. 
We define a subvaraiety $X$ of $W$ by the following equations: 
\begin{equation}
\begin{cases}
z^{3}-y \left( y^{2}\prod_{i=1}^{12}(u-a_{i})-x^{2} \right) =0, \\
z_{1}^{3}-y_{1} \left( y_{1}^{2}\prod_{i=1}^{12}(1-u_{1}a_{i})-x_{1}^{2} \right) =0
\end{cases} 
\tag{$\ast $}
\end{equation}
where $a_{i}$ ($i=1,2,\dots ,12$) are distinct complex numbers. 

Let $\pi $ be a projection from $X$ to the $u$-sphere $\mathbb{P}^{1}$. 
It is easy to see that  $X$ is a $K3$ surface and $\pi ^{-1}(u)$ is a non-singular elliptic curve with the 
functional invariant $0$ for every $u$ except $a_{i}$ ($i=1,\dots ,12$). 
Moreover we can see that $\pi ^{-1}(a_{i})$ is a singular fiber of type II. 
We define an automorphism $\varphi $ of $X$ as follows: $\varphi ([x:y:z], u )=([x:y:\zeta z], u )$ where 
$\zeta $ is a primitive third root of unity. Obviously $\varphi $ is of order $3$.
Now we remark that $X$ has a section $E$ defined by $y=0$.
Let $F$ be the class of a smooth fiber of $\pi $. 
Then $E$, $F+E$ generate the N\'{e}ron-Severi lattice of $X$ isometric to $U$.

The fixed locus of $\varphi $ is the set $\{z=0\}$. 
That is $\{y=0\}\coprod \{y^{2}\prod_{i=1}^{12}(u-a_{i})-x^{2}=0\}$. 
Clearly, the genus of the curve defined by $y=0$ is 0. 
Let $C$ be the curve defined by $y^{2}\prod_{i=1}^{12}(u-a_{i})-x^{2}=0$. 
The automorphism $\varphi $ induces an automorphism of order 3 on 
any smooth fiber of $\pi $ and $\pi ^{-1}(a_{i})$ ($i=1,\dots ,12$). 
Since $\varphi $ preserves a cusp of $\pi ^{-1}(a_{i})$, $C$ and $\pi ^{-1}(a_{i})$ intersect at the cusp of $\pi ^{-1}(a_{i})$. 
Thus we calculate the genus of $C$ by the Hurwitz formula: 
\[ 2g(C)-2 = 2(2g(\mathbb{P}^{1})-2)+12(2-1). \]
Therefore $X^{\varphi}=C^{(5)}\amalg \mathbb{P}^{1}$.

\vspace{5mm}

Next we consider that $\pi $ has no sections.

The next two examples are explained in detail by \cite{Kondo2}.

\begin{example}[\cite{Kondo2} $\S 2$] \label{exu3}
Let $C$ be a smooth non-hyperelliptic curve of genus 4. Then  its canonical model is the
complete intersection of an irreducible quadric surface $Q$ and an irreducible cubic surface $S$ in $\mathbb{P}^{3}$.
Let $X$ be the triple cover of $Q$ branched along $C$. 
Then $X$ is a $K3$ surface with an automorphism $\varphi $ of order $3$ and $C$ is a fix curve of $\varphi $.
Since $\varphi $ has a fixed curve $C$, $\varphi $ acts on $H^{0}(X, \Omega ^{2}_{X})$ 
as a multiplication by third root of unity (\cite{Nikulin2}, $\S 5$). 
Hence $\varphi $ is a non-symplectic automorphism.

Let $E$ (resp. $F$) be the inverse image of a smooth fiber of one of the rulings of $Q$ (resp. another ruling of $Q$). 
Then $E$, $F$ are elliptic curves with $\langle E, F \rangle =3$ and $E$, $F$ generate 
the N\'{e}ron-Severi lattice of $X$ isometric to $U(3)$.

Since $S_{X}=U(3)$ contains no $(-2)$ vectors, $X$ has no non-singular rational curves. 
And by the topological Lefschetz formula, $\chi _{\text{top}}(X^{\varphi})=-6=\chi _{\text{top}}(C)$.
Hence $X^{\varphi}$ contains no isolated points.
Thus the fixed locus of $\varphi $ is only $C$. 
Hence $X^{\varphi}=C^{(4)}$.
\end{example}

\begin{example}[\cite{Kondo2} $\S 3$]\label {node} 
Let $C_{1}$ be a curve in a smooth quadric $Q$ of bidegree $(3,3)$ with one node $p$. 
Let $L_{1}$, $L_{2}$ be the two lines through $p$. 
First blow up at $p$ and denote  by $E$ the exceptional curve. 
Next blow up the two points in which $E$ and the proper transform of $C_{1}$ meet. 
Then take the triple cover $X'$ branched along the proper transform of $C_{1}$ and $E$. 
Then $X'$ contains as exceptional curve of the first kind which is the pullback of the proper transform of $E$.
By contracting this exceptional curve to a smooth point $q$, we have a $K3$ surface $X_{1}$. 
Now for the generic $K3$ surface $X_{1}$, the N\'{e}ron-Severi lattice $S_{X_{1}}$ is isometric to $U(3)\oplus A_{2}$.

It is easy to see that the fixed locus of $\varphi $ is the smooth point $q$ and $\widetilde{C_{1}}$ where 
$\widetilde{C_{1}}$ is a smooth curve given by the normalization of $C_{1}$ at $p$. 
Since the genus of $\widetilde{C_{1}}=3$, we have $X_{1}^{\varphi}=\{q\}\amalg C^{(3)}$.
\end{example}

\begin{example}
Let $C_{2}$ be a curve in a smooth quadric $Q$ of bidegree $(3,3)$ with two nodes. 
By the same construction as  Example \ref{node}, we have a $K3$ surface $X_{2}$ with 
$S_{X_{2}}\simeq U(3)\oplus A_{2}^{\oplus 2}$ and a non-symplectic automorphism $\varphi $ of order 3.
And it is easy to see $X_{2}^{\varphi}=\{P_{1}\}\amalg \{P_{2}\}\amalg C^{(2)}$.

Similarly let $C_{3}$ be a curve in a smooth quadric $Q$ of bidegree $(3,3)$ with three nodes. 
Then we have a $K3$ surface $X_{3}$ with $S_{X_{3}}\simeq U(3)\oplus A_{2}^{\oplus 3}$ and 
a non-symplectic automorphism $\varphi $ of order 3.
And it is easy to see $X_{3}^{\varphi}=\{P_{1}\}\amalg \{P_{2}\}\amalg \{P_{3}\}\amalg C^{(1)}$.

Finally let $C_{4}$ be a curve in a smooth quadric $Q$ of bidegree $(3,3)$ with four nodes. 
Similarly we have a $K3$ surface $X_{4}$ with $S_{X_{4}}\simeq U(3)\oplus A_{2}^{\oplus 4}$ and 
a non-symplectic automorphism $\varphi $ of order 3.
Then we see $X_{4}^{\varphi}=\{P_{1}\}\amalg \{P_{2}\}\amalg \{P_{3}\}\amalg \{P_{4}\}\amalg C^{(0)}$.
\end{example}

\begin{example}[\cite{Zhang} Example 5.3.]\label{exz}
In this example, we construct a $K3$ surface $X$ with $S_{X}=U(3)\oplus E_{6}^{\ast }(3)$.

Denote by $[x : y : z ]$ the homogeneous coordinates of $\mathbb{P}^{2}$. 
Consider three cubic curves with a cusp of $\mathbb{P}^{2}$:
\[ C_{1}: x^{3}=y^{2}z, \ \ \ C_{2}:y^{3}=z^{2}x, \ \ \ C_{3}:z^{3}=x^{2}y. \]
Let $\eta$ be a primitive 7-th root of unity.
Then $C_{1}\cap C_{2}\cap C_{3}=\{[\eta ^{3j} :  \eta ^{j} : 1]|0\leq j\leq 6\}$.
Let $\pi:Y\longrightarrow \mathbb{P}^{2}$ be the blow-up at $[1 : 0 : 0]$, $[0 : 1 : 0]$, $[0 : 0 : 1]$, 
and seven points of $C_{1}\cap C_{2}\cap C_{3}$, let $D_{i}$ be the strict transform of $C_{i}$.
It is easy to see $D_{i}^{2}=-3$, $\chi_{\text{top}} (Y)=13$, and 
$0=\pi ^{\ast }\sum _{i=1}^{3}(C_{i}+3K_{\mathbb{P}^{2}})=\sum _{i=1}^{3}(D_{i}+3K_{Y})$.

Now let $\widetilde{Y}$ be the triple cover of $Y$ branched along the divisor $D_{1}+D_{2}+D_{3}$.
Note that $\widetilde{Y}$ has three ($-1$) curves. 
We contract these curves. 
Then we have a $K3$ surface $X$ with a non-symplectic automorphism of order $3$.
And the fixed locus consists of three isolated points. 
\end{example}

\end{document}